\documentclass[oneside,english]{amsart}
\usepackage[T1]{fontenc}
\usepackage[latin9]{inputenc}
\usepackage{amsthm}
\usepackage{amstext}
\usepackage{amssymb}
\usepackage{esint}
\usepackage{graphicx}
\usepackage{enumerate}
\usepackage{amscd}

\usepackage{graphicx}
\graphicspath{{noiseimages/}}

\makeatletter
\newtheorem{theorem}{Theorem}[section]

\numberwithin{figure}{section}
\theoremstyle{definition}

\newtheorem{example}[theorem]{Example}

\theoremstyle{remark}
\newtheorem{remark}[theorem]{Remark}

\numberwithin{equation}{section}
\makeatother

\usepackage{babel}

\begin{document}

\title[Neumann Eigenvalues]{Sobolev Extension Operators and Neumann Eigenvalues}

\author{V.~Gol'dshtein, V.~Pchelintsev, A.~Ukhlov}

\begin{abstract}
In this paper we apply estimates of the norms of Sobolev extension operators to the spectral estimates of of the first nontrivial Neumann eigenvalue of the Laplace operator in non-convex extension domains. As a consequence we obtain a connection between resonant frequencies of free membranes and the smallest-circle problem (initially proposed by J.~J.~Sylvester in 1857).
\end{abstract}
\maketitle

\section{Introduction}

Let $\Omega\subset\mathbb R^n$ be a bounded domain (an open connected set). Then $\Omega$ is called \textit{a Sobolev $L_2^1$-extension domain} if there exists a continuous linear extension operator 
\begin{equation}\label{CEO}
E:L_2^1(\Omega) \to L_2^1(\mathbb R^n)
\end{equation}
such that $E u|_{\Omega}=u$ for all $u \in L_2^1(\Omega)$.

Lipschitz domains are examples of Sobolev $L_2^1$-extension domains \cite{S70}. In \cite{GLV79} was given characterization of planar extension domains in terms of quasiconformal geometry, namely a simply connected planar domain $\Omega$ is a Sobolev $L_2^1$-extension domain if and only if $\Omega$ is an image of the unit disc $\mathbb D$ under a $K$-quasiconformal mapping of the plane $\mathbb R^2$ onto itself ($K$-quasidisc) for some $K\geq 1$. 

The aim of this article is estimate the Neumann eigenvalues of the Laplace operator defined in bounded domains $\Omega$ in terms of the norms of the extension operators. By the Min--Max Principle (see, for example, \cite{D95}) the first non-trivial Neumann eigenvalue $\mu_1(\Omega)$ for the Laplacian can be characterized as
$$
\mu_1(\Omega)=\min\left\{\frac{\int\limits _{\Omega}|\nabla u(x)|^{2}\, dx}{\int\limits _{\Omega}|u(x)|^{2}\, dx}:
u \in W^1_2(\Omega) \setminus \{0\},\,\, \int\limits _{\Omega}u\, dx=0 \right\}.
$$

Hence $\mu_1(\Omega)^{-\frac{1}{2}}$ is the best constant $B_{2,2}(\Omega)$ in the following Poincar\'e inequality
$$
\inf\limits _{c \in \mathbb R} \|u-c \mid L_2(\Omega)\| \leq B_{2,2}(\Omega) \|\nabla u \mid L_2(\Omega)\|, \quad u \in W_2^1(\Omega).
$$

It is well known (see, for example, \cite{M}) that if $\Omega\subset \mathbb R^n$ is the Sobolev $L_2^1$-extension domain then the spectrum of Neumann-Laplace operator in $\Omega$ is discrete and can be written in the form
of a non-decreasing sequence
\[
0=\mu_0({\Omega})<\mu_{1}({\Omega})\leq\mu_{2}({\Omega})\leq...\leq\mu_{n}({\Omega})\leq...\,.
\]

The main result of the article is

\vskip 0.3cm
\noindent
{\bf Theorem A.}
\textit{Let $\Omega\subset \mathbb R^n$ be a Sobolev $L_2^1$-extension domain. 
Then the following inequality holds 
\begin{equation}\label{Est}
\mu_1(\Omega) \geq \left(\frac{1}{\|E_{\Omega}\|}\frac{p_{n/2}}{R_{\Omega}}\right)^2 ,
\end{equation}
where $R_{\Omega}$ is a radius of a minimum enclosing ball $B_{\Omega}$ of $\Omega$, $p_{n/2}$ denotes the first positive zero of the function $\left(t^{1-n/2}J_{n/2}(t)\right)'$
and $\|E_{\Omega}\|$ denoted the norm of a continuous linear extension operator 
$$
E_{\Omega}: L^1_2(\Omega)\to L^1_2(B_{\Omega}).
$$
}

\begin{remark}
If a domain $\Omega$ has a center of a symmetry then $R_{\Omega}=d(\Omega)/2$, where $d(\Omega)=\sup\limits_{x,y\in\Omega}|x-y|$. In this case
\begin{equation}\label{EstSym}
\mu_1(\Omega) \geq \left(\frac{p_{n/2}}{d(\Omega)}\right)^2 \cdot \left(\frac{2}{\|E_{\Omega}\|}\right)^2.
\end{equation}

\end{remark}

The explicit value for the first non-trivial Neumann eigenvalue
of the Laplace operator in a $n$-ball of radius $R$ 
\begin{equation*}
\mu_1(B_R)=\left(\frac{p_{n/2}}{R}\right)^2,
\end{equation*}
where $p_{n/2}$ denotes the first positive zero of the function $\left(t^{1-n/2}J_{n/2}(t)\right)'$ (see, for example \cite{PS51, Weinb56}).
For some $n$, we give values of the parameter $p_{n/2}$ which represent first non-trivial Neumann eigenvalue in the unit ball in $\mathbb R^n$:
\vskip 0.3cm

\begin{center}
\begin{tabular}{|*{8}{p{20pt}|}}          \hline
\rule{0pt}{15pt}  $n$ & 2 & 3 & 4 & 5 & 6 & 7 & 8
\\          \hline
$p_{n/2}$  & 1.841 & 2.081 & 2.299 & 2.501 & 2.688 & 2,864 & 3.031        
\\          \hline
\end{tabular}
\end{center} 
\vskip 0.3cm

The proof of Theorem A is based on
\vskip 0.3cm
\noindent
{\bf Theorem B.} 
\textit{ Let $\Omega\subset \mathbb R^n$ be a Sobolev $L_2^1$-extension domain. Then for any bounded Lipschitz domain $\widetilde{\Omega}\supset \Omega$   
$$
\mu_1(\widetilde{\Omega}) \leq \|E_{\Omega}\|^2 \cdot \mu_1(\Omega),
$$
where $\|E_{\Omega}\|$ denoted the norm of a continuous linear extension operator 
$$
E_{\Omega}: L^1_2(\Omega)\to L^1_2(\widetilde{\Omega}).
$$
}


In the space case general estimates of the norm of the extension operator \eqref{CEO} is an open and complicated problem and we use estimates from \cite{Mikh} for balls and star-shaped domains. 

The construction of the extension operator which is based on Whitney decomposition \cite{Wh34} was studied by many authors
(see, for example, \cite{Cal61, Jon81,Shv16,Shv17,S70}), but in this case estimates of the norms of extension operators is an open complicated problem.

Extension operators on Sobolev spaces in terms of a measure density condition were studied in \cite{HKT08,HKT08_2}.

In the planar case we use the extension operator which was suggested in \cite{GLV79} and is based on quasiconformal reflections \cite{Ahl66} in Ahlfors domains ($K$-quasidiscs).
In this case 
\begin{equation*}
\|E_{\Omega}(u) \mid L^{1}_{2}(\widetilde{\Omega})\|\leq \|E(u) \mid L^{1}_{2}(\mathbb R^2)\| \leq (1+K)\|u \mid L^{1}_{2}(\Omega)\|,\,\,\text{for all}\,\,u\in L^{1}_{2}(\Omega).
\end{equation*}
Hence, $\|E_{\Omega}\|\leq 1+K$ and Theorem~A can be refined as

\vskip 0.3cm
\noindent
{\bf Corollary A.}
\textit{Let $\Omega$ be a $K$-quasidisc. Then 
\[
\mu_1(\Omega) \geq \left(\frac{j'_{1,1}}{R_{\Omega}}\right)^2 \cdot \left(\frac{1}{1+K}\right)^2,
\] 
where $R_{\Omega}$ is a radius of a minimum enclosing ball $B_{\Omega}$ of $\Omega$ and $j'_{1,1}\approx 1.84118$ denotes the first positive zero of the derivative of the Bessel function $J_1$.
}

\begin{remark}
Corollary A give a connection between resonant frequencies of free membranes and the smallest-circle problem. This problem was initially proposed by J.~J.~Sylvester in 1857. 
\end{remark}

The classical lower estimate for the first non-trivial eigenvalue
\begin{equation}
\label{eq:PW}
\mu_1(\Omega)\geq \frac{\pi^2}{d(\Omega)^2}
\end{equation}
was proved in \cite{PW} (see, also \cite{ENT, FNT}) for convex domains. 

\begin{remark}
If a domain $\Omega$ has a center of a symmetry and a norm of an extension operator 
$$
E_{\Omega}: L^1_2(\Omega)\to L^1_2(B_{\Omega})
$$
satisfies to the following condition 
$$
\|E_{\Omega}\|\leq \frac{2p_{n/2}}{\pi}
$$ 
then estimate \eqref{Est} is better than the classical estimate \eqref{eq:PW}.
\end{remark}

\begin{example}
Consider $n$-dimensional half-ball $B^{-}=\{x\in\mathbb R^n: |x|<1\,\,\&\,\,x_n<0\}$. Define the extension operator  
$$
E_{B^{-}}(u)=
\begin{cases}
u(x_1,...,x_{n-1},x_n), &\,\,\text{if}\,\,x\in B^{-},\\
u(x_1,...,x_{n-1},-x_n), &\,\,\text{if}\,\,x\in B^{+},
\end{cases}
$$
where $B^{+}=\{x\in\mathbb R^n: |x|<1\,\,\&\,\,x_n>0\}$.

Then $E_{B^{-}}:L^1_2(B^{-})\to L^1_2(B)$, $B=\{x\in\mathbb R^n: |x|<1\}$, and $\|E_{B^{-}}\|=\sqrt{2}$. Hence by Theorem A
$$
\mu_1(B^{-})\geq \frac{p^2_{n/2}}{2}>\frac{\pi^2}{4},\,\,\text{if}\,\,n\geq 4.
$$
So for $\mu_1(B^{-})$ the estimate by Theorem A is better than classical estimates~(\ref{eq:PW}) for $n\geq 4$.

\end{example}

In the present work we suggest lower estimates in non-convex domains in the terms of extension operators. In the planar case these estimates can be reformulated in the terms of (quasi)conformal geometry of $\Omega$. Another approach to the lower spectral estimates in non-convex domains is based on the geometric theory of composition operators on Sobolev spaces \cite{GPU17,GPU17_2,GU16,GU2016,GU17}.
 
Upper spectral estimates arise in works \cite{S54,Weinb56} and were intensively studied in the recent decades, see, for example, \cite{ A98, AB93, AL97, EP15, LM98}. 

\section{Extension operators}

Let $\Omega$ be a domain in $\mathbb{R}^n$, $n\geq2$.
For any $1\leq p<\infty$ we consider the Lebesgue space $L_p(\Omega)$ of measurable functions $u: \Omega \to \mathbb{R}$ equipped with the following norm:
\[
\|u\mid L_p(\Omega)\|=\biggr(\int\limits _{\Omega}|u(x)|^{p}\, dx\biggr)^{\frac{1}{p}}<\infty.
\]  

The Sobolev space $W^1_p(\Omega)$, $1\leq p<\infty$, is defined
as a Banach space of locally integrable weakly differentiable functions
$u:\Omega\to\mathbb{R}$ equipped with the following norm: 
\[
\|u\mid W^1_p(\Omega)\|=\biggr(\int\limits _{\Omega}|u(x)|^{p}\, dx\biggr)^{\frac{1}{p}}+
\biggr(\int\limits _{\Omega}|\nabla u(x)|^{p}\, dx\biggr)^{\frac{1}{p}}.
\]
Recall that the Sobolev space $W^1_p(\Omega)$ coincides with the closure of the space of smooth functions $C^{\infty}(\Omega)$ in the norm of $W^1_p(\Omega)$.

We consider also the homogeneous seminormed Sobolev space $L^1_p(\Omega)$, $1\leq p<\infty$,
of locally integrable weakly differentiable functions $u:\Omega\to\mathbb{R}$ equipped
with the following seminorm: 
\[
\|u\mid L^1_p(\Omega)\|=\biggr(\int\limits _{\Omega}|\nabla u(x)|^{p}\, dx\biggr)^{\frac{1}{p}}.
\]

\begin{remark}
By the standard definition functions of $L^1_p(\Omega)$ are defined only up to a set of measure zero, but they can be redefined quasieverywhere i.~e. up to a set of $p$-capacity zero. Indeed, every function $u\in L^1_p(\Omega)$ has a unique quasicontinuous representation $\tilde{u}\in L^1_p(\Omega)$. A function $\tilde{u}$ is called quasicontinuous if for any $\varepsilon >0$ there is an open  set $U_{\varepsilon}$ such that the $p$-capacity of $U_{\varepsilon}$ is less than $\varepsilon$ and on the set $\Omega\setminus U_{\varepsilon}$ the function  $\tilde{u}$ is continuous 
(see, for example \cite{HKM,M}). 
\end{remark} 

Recall that a continuous linear operator 
\[
E:L_p^1(\Omega) \to L_p^1(\mathbb R^n)
\]
satisfying the conditions
$$
Eu|_{\Omega}=u \quad \text{and} \quad \|E\| :=\sup\limits_{u \in L_p^1(\Omega)} \frac{\|Eu\mid{L_p^1(\mathbb R^n)}\|}{\|u\mid{L_p^1(\Omega)}\|}<\infty
$$
is called a continuous linear extension operator.

We say that $\Omega$ is \textit{a Sobolev $L_p^1$-extension domain} if there exists a continuous linear extension operator 
$$
E:L_p^1(\Omega) \to L_p^1(\mathbb R^n).
$$

It is well known that existence of an extension operator from $L_p^k(\Omega)$ to $L_p^k(\mathbb R^n)$, $k\in\mathbb N$, $n\geq2$, depends on the geometry of $\Omega$. In the case of Lipschitz domains, Calder\'on \cite{Cal61} constructed an extension operator on $L_p^k(\Omega)$ for $1 < p < \infty$, $k\in\mathbb N$. Stein \cite{S70} extended Calder\'on's result to include the endpoints $p=1, \infty$. 
Jones \cite{Jon81} introduced an extension operator on locally uniform domains. This is a much larger class of domains that includes examples with highly non-rectifiable boundaries. 

In \cite{GLV79} were obtained necessary and sufficient conditions for $L^1_2$-extensions from planar simply connected domains in terms of quasiconformal mappings and was obtained estimates of norms of extension operators.
Necessary and sufficient conditions for $L^1_p$-extension operators for $p>2$ from planar simply connected domains were obtained in \cite{Shv16, Shv17}  and were claimed  for $1<p<2$ in \cite{Kos98}. For case $p\neq 2$ estimates of norm of extension operators are unknown.

Using the extension operators theory on Sobolev spaces we estimates the weak domain monotonicity for the first non-trivial Neumann eigenvalues
in a bounded domains $\Omega \subset \mathbb R^n$.

The property of domain monotonicity holds for Dirichlet eigenvalues, i.e.,
if $\Omega \subset \widetilde{\Omega}$ are a bounded domains, then we have $\lambda_1(\Omega) \geq \lambda_1(\widetilde{\Omega})$.
This property does not holds for Neumann eigenvalues, even in the case of convex domains (see, for example, \cite{GN13}). For Sobolev $L_2^1$-extension domain there exists something like a quasi-monotonicity property that it can sees from the following theorem:

\vskip 0.3cm
\noindent
{\bf Theorem B.} 
\textit{ Let $\Omega\subset \mathbb R^n$ be a Sobolev $L_2^1$-extension domain. Then for any bounded Lipschitz domain $\widetilde{\Omega}\supset \Omega$   
$$
\mu_1(\widetilde{\Omega}) \leq \|E_{\Omega}\|^2 \cdot \mu_1(\Omega),
$$
where $\|E_{\Omega}\|$ denoted the norm of a continuous linear extension operator 
$$
E_{\Omega}: L^1_2(\Omega)\to L^1_2(\widetilde{\Omega}).
$$
}

\begin{proof}
Because $\Omega$ is the $L^1_2$-extension domain then 
there exists a continuous linear extension operator 
\begin{equation}\label{ExtOper}
E_{\Omega}:L_2^1(\Omega) \to L_2^1(\widetilde{\Omega})
\end{equation}
defined by formula
\[ (E_{\Omega}(u))(x) = \begin{cases}
u(x) & \text{if $x \in \Omega$,} \\
\widetilde{u}(x) & \text{if $x \in \widetilde{\Omega} \setminus {\Omega}$,}
\end{cases} \]
where $\widetilde{u}:\widetilde{\Omega} \setminus {\Omega} \to \mathbb R$ be an extension of the function $u$.

Hence for every function $u \in W_2^1(\Omega)$ we have
\begin{multline}\label{ineq1}
\|u-u_{\Omega} \mid L_2(\Omega)\| \\
{} =\inf\limits_{c\in \mathbb R}\|u-c\mid L_2(\Omega)\|=\inf\limits_{c\in \mathbb R}\|E_{\Omega}u-c\mid L_2(\Omega)\| \\
{} \leq \|E_{\Omega}u-(E_{\Omega}u)_{\widetilde{\Omega}}\mid L_2(\Omega)\| \leq \|E_{\Omega}u-(E_{\Omega}u)_{\widetilde{\Omega}} \mid L_2(\widetilde{\Omega})\|.
\end{multline}

Here $u_{\Omega}$ and $(E_{\Omega}u)_{\widetilde{\Omega}}$ are mean values of corresponding functions $u$ and $E_{\Omega}u$.

Because $\widetilde{\Omega}$ is a Lipschitz domain, then taking into account the classical Poincar\'e inequality in Lipschitz domains \cite{M}
$$
\inf\limits _{c \in \mathbb R} \|Eu-c \mid L_2(\widetilde{\Omega})\| \leq B_{2,2}(\widetilde{\Omega}) \|\nabla Eu \mid L_2(\widetilde{\Omega})\|, \quad u \in W_2^1(\widetilde{\Omega}), 
$$
and the continuity of the linear extension operator \eqref{ExtOper}, i.~e.,
$$
\|E_{\Omega}u \mid L_2^1(\widetilde{\Omega})\| \leq \|E_{\Omega}\| \cdot \|u \mid L_2^1(\Omega)\|, 
$$ 
we obtain
\begin{multline}\label{ineq2}
\|E_{\Omega}u-(E_{\Omega}u)_{\widetilde{\Omega}} \mid L_2(\widetilde{\Omega})\| \\
\leq B_{2,2}(\widetilde{\Omega}) \cdot \|\nabla (E_{\Omega}u) \mid L_2(\widetilde{\Omega})\| 
\leq B_{2,2}(\widetilde{\Omega}) \cdot \|E_{\Omega}\| \cdot \|\nabla u \mid L_2(\Omega)\|.
\end{multline}

Combining inequalities \eqref{ineq1} and \eqref{ineq2} we have
$$
\|u-u_{\Omega} \mid L_2(\Omega)\| \leq B_{2,2}(\Omega) \cdot \|\nabla u \mid L_2(\Omega)\|,
$$
where 
$$
B_{2,2}(\Omega) \leq B_{2,2}(\widetilde{\Omega}) \cdot \|E_{\Omega}\|.
$$

By the Min-Max Principle \cite{D95}  $\mu_1(\Omega)^{-1}=B_{2,2}^2(\Omega)$. Thus, we finally have
$$
\mu_1(\widetilde{\Omega}) \leq \|E_{\Omega}\|^2 \cdot \mu_1(\Omega).  
$$
\end{proof}

Let $B_{\Omega}$ be a minimum enclosing ball of $\Omega$. Taking in Theorem B $\widetilde{\Omega}=B_{\Omega}$  we obtain the main result of the article: 

\vskip 0.3cm
\noindent
{\bf Theorem A.} 
\textit{Let $\Omega\subset \mathbb R^n$ be a Sobolev $L_2^1$-extension domain. 
Then the following inequality holds 
\begin{equation*}
\mu_1(\Omega) \geq \left(\frac{1}{\|E_{\Omega}\|}\frac{p_{n/2}}{R_{\Omega}}\right)^2 ,
\end{equation*}
where $R_{\Omega}$ is a radius of a minimum enclosing ball $B_{\Omega}$ of $\Omega$, $p_{n/2}$ denotes the first positive zero of the function $\left(t^{1-n/2}J_{n/2}(t)\right)'$
and $\|E\|$ denoted the norm of continuous extension operator 
$$
E_{\Omega}: L^1_2(\Omega)\to L^1_2(B_{\Omega}).
$$
}


\subsection{Spectral estimates in planar domains}
In this section we give estimates of the norm of extension operators \eqref{ExtOper} in Ahlfors domains ($K$-quasidiscs).
Recall that a domain $\Omega\subset\mathbb R^2$ is called a $K$-quasidisc if it is an image of the unit disc $\mathbb D$ under a $K$-quasiconformal mapping of the plane $\mathbb R^2$ onto itself. Quasidiscs represent large class domains that  include some fractal type domains (like snowflakes).

Because there exists a M\"obius transformation of the unit disc $\mathbb D\subset\mathbb R^2$ onto the upper halfplane $H^{+}$ we can replace in the definition of quasidiscs $\mathbb D$ to $H^{+}$.
By the Riemann Theorem there exists a conformal mapping $\varphi$ of the upper halfplane $H^{+}$ onto $\Omega$.
Since $\Omega$ is a $K$-quasidisc there exists a quasiconformal extension of the conformal mapping $\varphi$ to a $K^2$-quasiconformal homeomorphism 
$\widetilde{\varphi} : \mathbb R^2 \to \mathbb R^2$ (see, for instance, \cite{Ahl66}).

By \cite{VG75} a homeomorphism $\psi : \Omega \to \widetilde{\Omega}$ is a $K$-quasiconformal mapping if and only if $\psi$ generates by the composition rule 
$\psi^{*}(\widetilde{u})= \widetilde{u} \circ \psi$ a bounded composition operator on Sobolev spaces $L_2^1(\Omega)$ and $L_2^1(\widetilde{\Omega})$:
$$
\|\psi^{*}(\widetilde{u}) \mid L_2^1(\Omega)\| \leq K^{\frac{1}{2}} \|\widetilde{u} \mid L_2^1(\widetilde{\Omega})\| 
$$
for any $\widetilde{u} \in L_2^1(\widetilde{\Omega})$.

Consider the following diagram
\[
\begin{CD}
L^{1}_{2}(\Omega) @>\varphi^{*}>> L^{1}_{2}(H^{+})\\
@V{(\varphi^{-1})^{*}\circ\omega\circ\varphi^{*}}VV  @VV{\omega}V \\
L^{1}_{2}(\mathbb R^2) @<(\varphi^{-1})^{*}<< L^{1}_{2}(\mathbb R^2)
\end{CD}
\]
where $\omega$ is a continuous extension operator (a symmetry with respect to the real axis)
that  extend any function $u \in L_2^1(H^{+})$ to a function $\widetilde{u}$ from $L_2^1(\mathbb R^2)$.

In accordance to this diagram we define the extension operator on Sobolev spaces
$$
E:L_2^1(\Omega) \to L^{1}_{2}(\mathbb R^2)
$$
by the formula
\[ (Eu)(x) = \begin{cases}
u(x) & \text{if $x \in \Omega$,} \\
\widetilde{u}(x) & \text{if $x \in \mathbb R^2 \setminus {\Omega}$,}
\end{cases} 
\]
where $\widetilde{u}:\mathbb R^2 \setminus {\Omega} \to \mathbb R$ is defined as $\widetilde{u}=(\varphi^{-1})^{*}\circ\omega\circ\varphi^{*}u$.

Now, using the above diagram, it is easy to check that 
\begin{multline*}
\|E u \mid L^{1}_{2}(\mathbb R^2)\| = \biggr(\int\limits _{\Omega}|\nabla u(x)|^{2}\, dx\biggr)^{\frac{1}{2}}
+ \biggr(\int\limits _{\mathbb R^2 \setminus \overline{\Omega}}|\nabla \widetilde{u}(x)|^{2}\, dx\biggr)^{\frac{1}{2}} \\
{} \leq (1+K)\biggr(\int\limits _{\Omega}|\nabla u(x)|^{2}\, dx\biggr)^{\frac{1}{2}}
=(1+K)\|u \mid L^{1}_{2}(\Omega)\|.
\end{multline*}

So, the norm of this extension operator $\|E\|\leq 1+K$.

Thus, taking into account Theorem A we have the following result:
\vskip 0.3cm
\noindent
{\bf Corollary A.}
\textit{Let $\Omega$ be a $K$-quasidisc. Then 
\[
\mu_1(\Omega) \geq \left(\frac{j'_{1,1}}{R_{\Omega}}\right)^2 \cdot \left(\frac{1}{1+K}\right)^2,
\] 
where $R_{\Omega}$ is a radius of a minimum enclosing ball $B_{\Omega}$ of $\Omega$ and $j'_{1,1}\approx 1.84118$ denotes the first positive zero of the derivative of the Bessel function $J_1$.
}
\vskip 0.3cm

As an example, we consider the following non-convex domain 
$$
\Omega := \left\{(x,y)\in \mathbb R^2 : -\frac{1}{2}<x<\frac{1}{2}, \quad \frac{1}{2}-|x|<y<\frac{3}{2}-|x|\right\}
$$
which is the image of the square
$$
Q:= \left\{(x,y)\in \mathbb R^2 : -\frac{1}{2}<x<\frac{1}{2}, \quad 0<y<1\right\}
$$
under a $K$-quasiconformal mapping of the plane $\mathbb R^2$ onto itself with the coefficient of quasiconformality $K=(3+\sqrt{5})/2$.

We show that domain $\Omega$ is a $K$-quasidisc. For this we input the following notation:
$$
Q=\left[-\frac{1}{2},\frac{1}{2}\right]\times[0,1], \quad S=\left[-\frac{1}{2},\frac{1}{2}\right]\times \mathbb R,
$$
$$
Q^{-}=\left[-\frac{1}{2},0\right]\times[0,1], \quad S^{-}=\left[-\frac{1}{2},0\right]\times \mathbb R,
$$
$$
Q^{+}=\left[0,\frac{1}{2}\right]\times[0,1], \quad S^{+}=\left[0,\frac{1}{2}\right]\times \mathbb R. 
$$
Let mapping $\widetilde{\varphi}_{+}:\mathbb R^2 \to \mathbb R^2$ is defined by the rule $\widetilde{\varphi}_{+}:=\left\{x, x+y+\frac{1}{2}\right\}$ 
and mapping $\widetilde{\varphi}_{-}:\mathbb R^2 \to \mathbb R^2$ is defined by the rule $\widetilde{\varphi}_{-}:=\left\{x, y-x+\frac{1}{2}\right\}$.
It is easy to see that
$$
\widetilde{\varphi}_{+}(0,y)=\left\{0,y+\frac{1}{2}\right\}=\widetilde{\varphi}_{-}(0,y)
$$
and
$$
\widetilde{\varphi}_{+}\left(-\frac{1}{2},y\right)=\left\{-\frac{1}{2},y\right\}, \quad 
\widetilde{\varphi}_{-}\left(\frac{1}{2},y\right)=\left\{\frac{1}{2},y\right\}.
$$
Let mapping $\varphi:Q\to \mathbb R^2$ is defined as
$$
\varphi|_{Q^{-}}:=\varphi_{+} \quad \text{and} \quad \varphi|_{Q^{+}}:=\varphi_{-},
$$
and mapping $\overline{\varphi}: S \to \mathbb R^2$ is defined as
$$
\overline{\varphi}|_{S^{-}}:=\overline{\varphi}_{+} \quad \text{and} \quad \overline{\varphi}|_{S^{+}}:=\overline{\varphi}_{-},
$$
where
$$
\varphi_{+}:=\widetilde{\varphi}_{+}|_{Q^{-}}, \quad \overline{\varphi}_{+}:=\widetilde{\varphi}_{+}|_{S^{-}}, \quad
\varphi_{-}:=\widetilde{\varphi}_{-}|_{Q^{+}}, \quad \overline{\varphi}_{-}:=\widetilde{\varphi}_{-}|_{S^{+}}.
$$
Then extended mapping $\widetilde{\varphi}:\mathbb R^2 \to \mathbb R^2$ is defined as
$$
\widetilde{\varphi}|_{Q}:=\varphi, \quad \widetilde{\varphi}|_{S}:=\overline{\varphi}, \quad
\text{and} \quad \widetilde{\varphi}|_{R^2 \setminus S}:= \text{id}.
$$

Now we calculate the quasiconformality coefficient for the mapping $\varphi:Q\to \mathbb R^2$,
using the following formula
$$
K=\frac{\lambda}{J_{\varphi}(x,y)}.
$$
Here $\lambda$ is the largest eigenvalue of the matrix $A=DD^T$, where 
$D=D\varphi(x,y)$ is Jacobi matrix of mapping $\varphi=\varphi(x,y)$ and $J_{\varphi}(x,y)=\det D\varphi(x,y)$ is its Jacobian.

Note that the Jacobi matrix corresponding to the mapping $\varphi=\varphi(x,y)$ on $S^{-}$, $S^{+}$ and outside of $S$, respectively, has the form
$$
D=\begin{pmatrix} 1 & 0 \\ 1 & 1 \end{pmatrix}, \quad
D=\begin{pmatrix} 1 & 0 \\ -1 & 1 \end{pmatrix}, \quad \text{and} \quad
D=\begin{pmatrix} 1 & 0 \\ 0 & 1 \end{pmatrix}.
$$

A straightforward calculation yields
$$
J_{\varphi}(x,y)=1 \quad \text{and} \quad \lambda=\frac{3+\sqrt{5}}{2}.
$$
Hence
$$
K=\frac{\lambda}{J_{\varphi}(x,y)}=\frac{3+\sqrt{5}}{2}.
$$

Thus, we established that the domain $\Omega$ is $K$-quasidisc with the quasiconformality coefficient $K=(3+\sqrt{5})/2$.

It is easy to show that
$$
d(\Omega)=\sup\limits_{x,y\in\Omega}|x-y|=\frac{\sqrt{10}}{2}.
$$

Then by Corollary~A we have
\[
\mu_1(\Omega) \geq \left(\frac{j'_{1,1}}{d(\Omega)}\right)^2 \cdot \left(\frac{2}{1+K}\right)^2 \approx \frac{2}{5}.
\] 

\vskip 0.3cm

As another applications of Corollary~A, we obtain the lower estimates of the first non-trivial eigenvalue on the 
Neumann eigenvalue problem for the Laplace operator in the star-shaped and spiral-shaped domains.

\vskip 0.3cm

\textbf{Star-shaped domains.} We say that a domain $\Omega^*$ is $\beta$-star-shaped (with respect to $z_0=0$)
if the function $\varphi(z)$, $\varphi(0)=0$, con formally maps a unit disc $\mathbb{D}$ onto $\Omega^*$ and the condition satisfies \cite{FKZ}:
\[
\left|\arg \frac{z \varphi^{\prime}(z)}{\varphi(z)} \right| \leq \beta \pi/2, \quad 0 \leq \beta <1, \quad |z|<1.
\]

In \cite{FKZ} proved the following: 
the boundary of the $\beta$-star-shaped domain $\Omega^*$ is a $K$-quasicircle with $K=\cot ^2(1-\beta)\pi/4$.

Then by Corollary~A we have
\[
\mu_1(\Omega^*) \geq 4 \sin^4\left((1-\beta)\pi/4\right) \cdot \left(\frac{j'_{1,1}}{d(\Omega^*)}\right)^2.
\]

For example, the diffeomorphism 
$$
\varphi(z)=\tan z, \quad z=x+i y, \quad |z|<1, 
$$ 
conformally maps a unit disc onto $\frac{1}{2}$-star-shaped domain $\Omega^*$:

\begin{figure}[h!]
\centering
\includegraphics[width=0.4\textwidth]{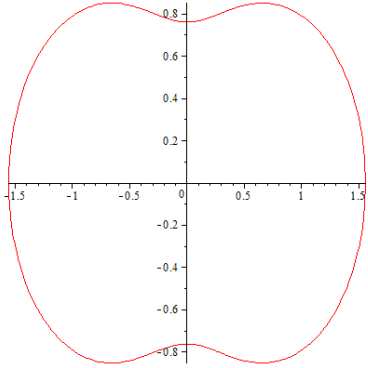}
\caption{Image of $\mathbb{D}$ under $\varphi(z)=\tan z$.}
\end{figure}

By straightforward calculation we get
\[
\mu_1(\Omega^*) \geq 4 \sin^4(\pi/8) \cdot \left(\frac{j'_{1,1}}{d(\Omega^*)}\right)^2
\approx 4 \left(\frac{\sqrt{2 - \sqrt{2}}}{2}\right)^4 \cdot \left(\frac{1.84118}{3.2}\right)^2 \approx \frac{1}{5}.
\]

\vskip 0.3cm

\textbf{Spiral-shaped domains.} We say that a domain $\Omega_s$ is $\beta$-spiral-shaped (with respect to $z_0=0$)
if the function $\varphi(z)$, $\varphi(0)=0$, conformally maps a unit disc $\mathbb{D}$ onto $\Omega_s$ and the condition satisfies \cite{S87, Sug}:
\[
\left|\arg e^{i \gamma} \frac{z \varphi^{\prime}(z)}{\varphi(z)} \right| \leq \beta \pi/2, \quad 0 \leq \beta <1, \quad |\gamma|<\beta \pi/2, \quad |z|<1.
\]

In \cite{S87, Sug} proved the following: the boundary of the $\beta$-spiral-shaped domain $\Omega_s$ is a $K$-quasicircle with $K=\cot ^2(1-\beta)\pi/4$.

Then by Corollary~A we have
\[
\mu_1(\Omega_s) \geq 4 \sin^4\left((1-\beta)\pi/4\right) \cdot \left(\frac{j'_{1,1}}{d(\Omega_s)}\right)^2.
\] 

\vskip 0.3cm

\subsection{Spectral estimates in space domains}
At the best of our knowledge estimates of norms of extension operators for $L^1_2$  in space domains are not known.

We found only Mikhlin's estimates \cite{Mikh} for balls and star-shaped domains.
Let $B_1$ and $B_R$ be balls in Euclidean space $\mathbb R^n$ with a common center at the origin 
and with radii $1$ and $R$, $R>1$, respectively. In this case, Mikhlin~\cite{Mikh} established that the norm of extension operator
\begin{equation}\label{EOM}
E_R:W_2^1(B_1) \to W_2^1(B_R)
\end{equation}
is 
\begin{equation}\label{EOE}
\|E_R\|^2=1+ \frac{I_{\alpha}(1)}{I_{\alpha+1}(1)} \cdot 
\frac{I_{\alpha}(R) K_{\alpha+1}(1)+K_{\alpha}(R) I_{\alpha+1}(1)}{I_{\alpha}(R) K_{\alpha}(1)-K_{\alpha}(R) I_{\alpha}(1)}.
\end{equation}
Here $\alpha=(n-2)/2$, $n>2$ and $I_{\nu}(z)$, $K_{\nu}(z)$ are the Bessel functions of the imaginary argument:
$$
I_{\nu}(z)=\sum_{m=0}^{\infty} \frac{1}{m! \Gamma(m+\nu+1)}\left(\frac{z}{2}\right)^{2m+\nu},
$$
where $\Gamma(\cdot)$ is the gamma function and
$$
K_{\nu}(z)=\frac{\pi}{2} \cdot \frac{I_{-\nu}(z)-I_{\nu}(z)}{\sin(\nu \pi)}.
$$

For $n=3$ and $R=2$ or $R=3$ we have the following values of the norm of extension operator \eqref{EOM}:
$$
\|E_2\|^2\approx8.38905 \quad and \quad \|E_3\|^2\approx7.50825.
$$

Let $u(x)$ be a nonnegative continuous function of class $W_{\infty}^1$ in some two-sided 
neighborhood of the unit sphere $S:=\left\{x:|x|=1\right\}$.
We set $\rho=|x|$, $\theta=x/\rho$ and denote
$$
M_1=\min\limits_{x \in S} u(x)>0, \quad M_2=\max\limits_{x \in S} u(x), \quad M_3=\sup\limits_{x \in S}|\nabla u(x)|.
$$

Let $\Omega_1$ and $\Omega_R$, $R>1$, be a star-shaped domains having the form
$$
\Omega_{\beta}:=\left\{x: \rho< \beta u(\theta)\right\}, \quad \beta >0.
$$ 
Then in \cite{Mikh} obtained the estimate of the norm of extension operator
\begin{equation}\label{EOM1}
E^*:W_2^1(\Omega_1) \to W_2^1(\Omega_R)
\end{equation}
which have the form
\begin{equation}\label{EOM2}
\|E^*\|^2 \leq 1+\left(\frac{M_2}{M_1}\right)^2\left(N_1 N_2\right)^2\left(\|E_R\|^2-1\right),
\end{equation}
where
$$
N_1^2=\max\left\{\frac{M_1^2+(n-1)M_3^2}{M_1^4},\,\frac{2}{M_1^2},\,1\right\}
$$
and
$$
N_2^2=\max\left\{M_2^2+2(n-1)M_3^2,\,2M_2^2,\,1\right\}.
$$

Hence, by Theorem~A  we have the following quasi-monotonicity result for Neumann eigenvalues in star-shaped domains:
\vskip 0.3cm
\noindent
{\bf Corollary B.}
\textit{Let $\Omega_1$ and $\Omega_R$, $R>1$, be star-shaped domains in Euclidean space $\mathbb R^n$.
Then the following quasi-monotonicity estimate  
\[
\mu_1(\Omega_R) \leq \mu_1(\Omega_1) \|E^*\|^2 \leq 1+\left(\frac{M_2}{M_1}\right)^2\left(N_1 N_2\right)^2\left(\|E_R\|^2-1\right)
\] 
holds.
}

\vskip 0.3cm

{\bf Acknowledgment.} 
The financial support by United States-Israel  Binational Science foundation (Grant 2014055) is gratefully acknowledged.

\vskip 0.3cm

Department of Mathematics, Ben-Gurion University of the Negev, P.O.Box 653, Beer Sheva, 8410501, Israel 
 
\emph{E-mail address:} \email{vladimir@math.bgu.ac.il} \\           
       
 Department of Higher Mathematics and Mathematical Physics, Tomsk Polytechnic University, 634050 Tomsk, Lenin Ave. 30, Russia;
 Department of General Mathematics, Tomsk State University, 634050 Tomsk, Lenin Ave. 36, Russia

 \emph{Current address:} Department of Mathematics, Ben-Gurion University of the Negev, P.O.Box 653, 
  Beer Sheva, 8410501, Israel  
							
 \emph{E-mail address:} \email{vpchelintsev@vtomske.ru}   \\
			  
	Department of Mathematics, Ben-Gurion University of the Negev, P.O.Box 653, Beer Sheva, 8410501, Israel 
							
	\emph{E-mail address:} \email{ukhlov@math.bgu.ac.il  

\end{document}